\documentclass{article}

\usepackage{wollatex}
\usepackage{format-short}

\newcommand{\gen}{\mathrm{gen}}
\newcommand{\Vt}{\widetilde{V}}
\newcommand{\ellt}{\tilde{\ell}}
\newcommand{\Oom}{\mathcal{O}}

\renewcommand{\P}{\mathcal{P}}
\newcommand{\Pf}{\P(\future)}
\newcommand{\PPf}{\P\(\P(\future)\)}
\newcommand{\Pinv}{\P_\mathrm{s}}
\newcommand{\Pe}{\P_\mathrm{e}(\D^\Z)}
\newcommand{\PO}{\Pinv(\D^\Z)}
\newcommand{\Ms}{\mathcal{M}}
\newcommand{\Mp}{\Ms_+}
\newcommand{\Msf}{\Ms(\future)}
\newcommand{\MCf}{\Ms_\C(\future)}

\newcommand{\Ph}{\widehat{P}}
\newcommand{\D}{\Delta}
\newcommand{\future}{\D^\N}
\newcommand{\A}{\mathcal{A}}
\newcommand{\shift}{\sigma}

\newcommand{\tauD}{\tau^\D}
\newcommand{\tauA}{\tau^\A}

\newcommand{\caus}{\mathfrak{C}}

\newcommand{\mueps}{\mu_\caus^P}
\newcommand{\Ceps}{C_\caus}	

\newcommand{\causup}{\mathfrak{C}}
\newcommand{\causupp}{\causup_P}

\newcommand{\Xp}[1][X]{#1_{-\N_0}}
\newcommand{\Xf}[1][X]{#1_{\N}}

\newcommand{\xp}[1][x]{\Xp[#1]}

\newcommand{\unita}{\unit_\A}
\newcommand{\unitaf}{\unit_{\Af}}
\newcommand{\tensor}[2]{#1^{\otimes #2}}
\newcommand{\AI}[1][I]{\A_{#1}}
\newcommand{\AZ}{\AI[\Z]}
\newcommand{\Af}{\AI[\N]}
\newcommand{\Ap}{\AI[-\N_0]}
\newcommand{\Afin}{\A_{<\mspace{-1mu}\infty}}
\newcommand{\State}{\mathcal{S}}
\newcommand{\Sf}{\State(\Af)}
\newcommand{\Se}{\State_\mathrm{e}}
\newcommand{\Ss}{\State_\mathrm{s}}
\newcommand{\Sef}{\Se(\Af)}
\newcommand{\Ssf}{\Ss(\Af)}
\newcommand{\B}{\mathfrak{B}}
\newcommand{\cont}{\mathcal{C}}

\renewcommand{\H}{\mathcal{H}}
\newcommand{\Hphi}{\H_\vphi}
\newcommand{\piphi}{\pi_\vphi}
\newcommand{\xiphi}{\xi_\vphi}

\newcommand{\bscalphi}[2]{\bscalprod{#1}{#2}_\vphi}

\newcommand{\rep}{\iota_\vphi}

\title{Process Dimension of Classical and Non-Commutative Processes}
\author{Wolfgang L\"ohr$^{1,2}$ \and Arleta Szko{\l}a$^1$ \and Nihat Ay$^{1,3}$}
\date{\today}

\begin{document}

\footnotetext[1]{Max Planck Institute for Mathematics in the Sciences, Inselstra{\ss}e 22, 04103 Leipzig, Germany}
\footnotetext[2]{Universit\"at Duisburg-Essen, Universit\"atsstra{\ss}e 2, 45141 Essen, Germany}
\footnotetext[3]{Santa Fe Institute, 1399 Hyde Park Road, Santa Fe, New Mexico 87501, USA}

\maketitle

\begin{abstract}
	We treat observable operator models (OOM) and their
	non-commutative generalisation, which we call NC-OOMs. A
	natural characteristic of a stochastic process in the context
	of classical OOM theory is the process dimension. We
	investigate its properties within the more general
	formulation, which allows to consider process dimension as a
	measure of complexity of non-commutative processes: We prove
	lower semi-continuity, and derive an ergodic decomposition
	formula. Further, we obtain results on the close relationship
	between the canonical OOM and the concept of causal states
	which underlies the definition of statistical complexity. In
	particular, the topological statistical complexity, i.e.\ the
	logarithm of the number of causal states, turns out to be an
	upper bound to the logarithm of process dimension.
\end{abstract}

\section{Introduction}

The main idea behind various complexity measures, such as statistical
complexity, is the same that gave rise to the famous Kolmogorov complexity.
Namely, the complexity is the ``size'' of some minimal ``representation'' of the
object of interest.  Different complexity measures are based on different exact
definitions of these terms. For Kolmogorov complexity, for instance,
representations are Turing machine programs computing individual binary strings,
and the size is their length. For statistical complexity (\cite{epsdef}), on the
contrary, the objects of interest are probability distributions of stochastic
processes instead of individual strings, and the representations are particular
kinds of predictive models in the sense of partially deterministic hidden Markov
models (HMM). Their size is measured by the Shannon entropy of the internal
states of the model.

Observable operator models (OOM) are generative algebraic models that represent
a stochastic process. The natural measure of size of an OOM is the dimension of
the corresponding real vector space. It is minimal for canonical OOMs of a given
stochastic process and was already identified as a characteristic of the process
called \emph{process dimension}. In the present contribution, following the
above mentioned reasoning, we want to consider the process dimension as a
complexity measure for stochastic processes. We give further indication that
this might be appropriate.
First, we show the close relation of the canonical OOMs to the concept of causal
states which are used to define statistical complexity. Second, we prove that
the process dimension, considered as function of the process, is lower
semi-continuous. Although there exists no generally accepted axiomatic
characterisation of functionals on the space of stochastic processes that
quantify complexity, we argue that every complexity measure should feature this
property.  Indeed, it would be strange to consider a process complex if there is
an approximating sequence with (uniformly) simple processes.
The natural topology for processes is in this context {\weaks} topology as
opposed to the much stronger variational topology, and lower semi-continuity
w.r.t.\ {\weaks} topology is a much stronger result.

The construction of causal states relies heavily on conditional probabilities.
This makes it difficult to extend the corresponding notion of statistical
complexity to the domain of non-commutative processes understood to be states on
a quasi-local C*-algebra. The algebraic formulation of OOMs, however, allows to
extend the concept of process dimension to the non-commutative setting. Indeed,
the construction of \emph{finitely correlated states} introduced by Fannes et
al.\ in \cite{Fannes} provides OOMs for a class of shift-invariant states on a
quasi-local C*-algebra. In the literature, these states are also known as
\emph{algebraic states}. In this paper, however, we refer to their original
name. We show lower semi-continuity of the process dimension also in this more
general setting.

The outline of our paper is as follows. In Section~2 we present main concepts in
the context of classical OOM theory using the terminology introduced by Herbert
Jaeger.  In particular, process dimension is defined.
In Section~3 we review the concept of statistical complexity and the underlying
notion of causal states which are defined for classical stationary process. In
\propref{closeq} we specify the way in which process dimension and causal states
are related. As direct implication we obtain an upper bound for, and an ergodic
decomposition of the process dimension in the classical case
(\correfs{upperbound}{classergdecomp}).
In Section 4 and 5 we treat non-commutative extensions. We start with the
definition of NC-OOMs referring to finitely correlated states. As in the
classical special case the corresponding process dimension is naturally
associated with canonical NC-OOMs. Our main results are contained in Section~5.
There we prove lower-semicontinuity and provide an ergodic decomposition formula
for the process dimension in the general case.

\section{Classical OOMs (Stochastic Modules)}

Fix a \emph{finite} set $\Delta$. We consider $\D$-valued stochastic processes $X_\N \defeq (X_k)_{k\iN}$, described by their
distributions $P\in\Pf$, and stationary processes $X_\Z$, described by their shift-invariant distributions $P\in\PO$. 

In \cite{Heller}, Alex Heller introduced a generalisation of functions of Markov chains, called \emph{stochastic modules}.
Later, Herbert Jaeger extended and reformulated this theory in the language of linear algebra (\cite{Jaeger}). We use his
terminology.

\begin{definition}\deflabel{oom}
	An \define{observable operator model} (\define{OOM}) with alphabet $\Delta$ is a quadruple $(V, T, v, \ell)$, where $V$
	is a real vector space, $T\colon \D\times V\to V$ is linear in the second argument, $v$ is an element of $V$, and $\ell$
	is a linear form on $V$, such that for $T_d(v) \defeq T(d, v)$, $n\iN$ and $d_1,\ldots,d_n \in\Delta$,
	\begin{gather*}
		\eqnitem{1}	\ell(v) \= 1,
		\sepeqnitem{2}	\ell\circ\sum_{d\in\Delta}  T_d \= \ell,
		\sepeqnitem{3}	P_{d_1,\ldots,d_n}\defeq \ell\circ T_{d_n}\circ\cdots\circ T_{d_1}(v) \ges 0.
	\end{gather*}
	The vector $v$ is called \define{initial vector}, the operators $T_d$ are called \define{observable operators} and
	the linear form $\ell$ is called \define{evaluation form}\footnote{Jaeger fixes a basis of $V$ instead of an evaluation
	form and defines $\ell$ to be the sum of coefficients in the basis expansion.}.
	The process $P\in\Pf$, defined by $P\([d_1,\cdots, d_n]\) \defeq P_{d_1,\ldots, d_n}$ ($n\iN$), is called
	\defemph{generated} by the OOM and the dimension $\dim(V)$ of $V$ is called \define{dimension} of the OOM. 
\end{definition}

It is easy to check that the $P_{d_1,\ldots, d_n}$ are a consistent set of finite-dimensional probabilities.
Therefore, by the Kolmogorov extension theorem, the process $P$ well-defined.
Every hidden Markov model (HMM)
with $n$ internal states canonically induces an $n$-dimensional OOM. For more details, see \cite{Jaeger}.

\begin{remark}
	More generally, an HMM with set $\Gamma$ of internal states can be interpreted as OOM with vector space $\Ms(\Gamma)$
	of signed measures of bounded variation on $\Gamma$.
\end{remark}

There is a canonical construction of an OOM of a given process $P \in \Pf$.
Let $\Msf$ be the space of signed measures of bounded variation on $\future$,
i.e.\ $\Msf = \linhull\(\Pf\)$, where $\linhull$ denotes the linear hull. Define the linear maps $\tauD_d\colon \Msf
\to \Msf$ by
	\[ \tauD_d(\mu) \defeq \mu\([d]\cap \shift^{-1}(\fdot)\) \]
where $\shift$ is the left-shift on $\future$. Further define $\ell_\D\colon \Msf \to \R$ by $\ell_\D(\mu)=\mu(\future)$, i.e.\ the
evaluation form $\ell_\D$ associates to a measure its total mass. For convenience we define
	\[ \tauD_{d_1\cdots d_n}\defeq \tauD_{d_n}\circ\cdots\circ \tauD_{d_1}. \]

\begin{definition}\deflabel{canoom}
	For $P\in\Pf$, let
		\[ Q_P \defeq \bset{\tauD_{d_1\cdots d_n}(P)}{n\iN_0,\;d_1,\ldots,d_n\in\D}  \und V_P \defeq \linhull(Q_P). \]
	For $d\in\D$, denote the function $V_P\to V_P$, $\mu\mapsto\tauD_d(\mu)$ with a slight abuse of notation again by
	$\tauD_d$.  Set $\tauD(d,\mu)\defeq \tauD_d(\mu)$.
	Then $(V_P, \tauD, P, \ell_\D)$ is called \define{canonical OOM} of $P$.
\end{definition}

Since $\tauD_d(V_P) \subseteqs V_P$, the canonical OOM is well-defined, and it generates $P$. It
has minimal dimension among all OOMs generating $P$, and is, up to isomorphism, unique with this minimality
(see \cite{Jaeger}). In particular, the dimension of $V_P$ is not bigger (but may be essentially smaller)
than the minimal number of internal states required for any HMM generating $P$.
Another characterisation of $V_P$ is in terms of conditional probabilities:
\[ V_P \= \linhull\, \Bset{\bcondp{\shift^{-n}(\fdot)}{[d_1,\ldots,d_n]}}
		{n\iN_0,\; d_1,\ldots,d_n\in\D,\;P\([d_1,\ldots,d_n]\)>0} . \]
This is true because if we normalise $\tauD_d$ pointwise, we obtain the corresponding conditional probability.
Note that the normalised version of $\tauD_d$ is not linear.

If $A$ is a finite dimensional cylinder set, the same holds for $[d]\cap \shift^{-1}(A)$. Therefore, $\tauD_d$
is {\weaks} continuous and, consequently, $\tauD_d$ maps the \weaks\ closure $\wsclosure{V_P}$ to itself.
Thus $\(\wsclosure{V_P}, \tauD, P, \ell_\D\)$ is an OOM of $P$, which we call the \define{closed canonical OOM}.
In the case of finite process dimension, which we are mostly interested in, the canonical OOM and the closed
canonical OOM coincide. For considering infinite pasts in the following section, however, the {\weaks} closure
of $V_P$ in $\Msf$ plays a crucial role.

\begin{example}[Canonical OOM and Shift HMM]
	The (one-sided) shift HMM of $P\in\Pf$ is a deterministic HMM with set $\Gamma\defeq\future$ of internal states. It
	is in general by no means minimal and it is not possible to restrict it to a smaller subset of $\Gamma$ such that it
	still generates $P$. If we interpret it as OOM, the internal Vector space is $V=\Msf$, the associated operators $T_d$
	are equal to the canonical ones, i.e.\ $T_d=\tauD_d$, the initial vector is the initial distribution of the shift
	HMM, i.e.  $v=P$ and the evaluation form is $\mu\mapsto\mu(\future)$. Now it is obvious that we can reduce every OOM
	to a ``cyclic'' version by restricting $V$ to $\linhull\bset{T_{d_1\ldots d_n}(v)}{n\iN_0,\; d_1,\ldots,d_n\in\D}$.
	This reduced shift OOM is just the canonical OOM and thus minimal, but it can in general not be interpreted as HMM.
\end{example}

\begin{definition}\deflabel{procdim}
	The \define{process dimension} of $P\in\Pf$ is the dimension of its canonical OOM $(V_P, \tauD, P, \ell_\D)$:
		\[ \dim(P) \defeq \dim(V_P)\ins\N\cup\sset{\infty}. \]
\end{definition}

The process dimension is derived from a canonical construction and at the same time the minimal dimension necessary for an
OOM-representation. Therefore, it is an important internal characteristic of the process and might be considered a complexity
measure. As we see in the following section, it is related to statistical complexity.

\section{Infinite Pasts and Causal States}

Now assume we are interested in a stationary process with infinite past, i.e.\ $X_\Z$, described by its
shift-invariant distribution $P\in\PO$. For simplicity assume that $X_k$ is the canonical projection on
$\D^\Z$. We define the canonical OOM of $P$ to be the canonical OOM of its restriction to positive times, i.e.
if $P_\N \defeq P\circ X_\N^{-1}$ is the distribution of $X_\N$ then
	\[ \dim(P) \defeq \dim(P_\N) \= \dim(V_{P_\N}). \]
In \cite{epsdef}, the causal states of such a stationary process were introduced and used to define
statistical complexity. Causal states are equivalence classes of past trajectories, where two of them are
identified if they induce the same conditional probability distribution on the future $X_{\N}$.
In this paper we prefer the alternative definition of causal states, where they are measures on the future, i.e.\
elements of $\Pf$. This viewpoint was introduced in \cite{funcpaper}.

If we observe the past $\Xp$ of $X_\Z$, the observation $\Xp=\xp$ induces a certain conditional probability
distribution $\condp{\Xf}{\Xp=\xp}\in\Pf$ on the future $\Xf$ of the process.  The causal state distribution
of $P\in\PO$ is the distribution of these conditional probabilities. In particular, it is a measure on
measures. More precisely we define

\begin{definition}
	Let $P\in\PO$. The \define{causal state distribution} $\mueps\in\PPf$ of $P$ is defined by 
		\[ \mueps \defeq P \circ \(\condp{\Xf}{\Xp}\)^{-1}, \]
	where we consider $\condp{\Xf}{\Xp}$ to be a (measurable) function from $\D^\Z$ to $\Pf$. Further, we define
		\[ \causupp \defeq \supp\(\mueps\) \subseteqs \Pf. \]
\end{definition}

\begin{enremark}
	\item $\mueps$ is the distribution of the $\Pf$-valued random variable $\condp{\Xf}{\Xp}$.
	\item The causal states correspond to the elements in the image of $\condp{\Xf}{\Xp}$. Thus, the set of causal states
		depends on the version of conditional probability. $\causupp$, on the other hand, is independent of the choice
		of conditional probability.
	\item The statistical complexity $\Ceps(P) \defeq H(\mueps)$ is the (Shannon) entropy of the causal state
		distribution. It was originally introduced by Grassberger as \emph{true measure complexity} in
		\cite{grassberger}.
\end{enremark}

We obtain the following relation between the causal state distribution and the closed canonical OOM vector space.
The {\weaks} closure of the canonical OOM vector space is equal to the {\weaks} closure of the vector space spanned
by the support $\causupp$ of the causal state distribution. In the finite dimensional case, this means that the
two vector spaces are equal. Because the OOM vector space is defined with finite-length pasts and infinite pasts
are used for the definition of the causal state distribution, we can interpret this result as follows. Unlike the
set of causal states, the canonical OOM vector space is the same if we consider finite or infinite pasts, provided
it is finite dimensional. For the infinite dimensional case, the situation is more subtle (see
\exref{infinitedim}). Note that $\wsclosure{V}$ denotes the closure of $V$ w.r.t.\ the {\weaks} topology. Recall
that $V_P=\linhull(Q_P)$.

\begin{proposition}\proplabel{closeq}
	Let $P\in\PO$. Then
		\[ \wsclosure{V_P} \= \wsclosure{\linhull(\causupp)}. \]
	In particular, because finite-dimensional spaces are closed, $\dim(P)\=\dim\(\linhull(\causupp)\)$.
\end{proposition}
\begin{proof}
\wcase{``$\subseteq$''} Let $\nu\in Q_P$. Then $\nu=\tauD_{d_1\cdots d_n}(P_\N)$ for some $d_1,\ldots,d_n\in\D$.
	Define $A$ to be the event that the past is $d_1,\ldots, d_n$, i.e.\ $A\defeq\sset{X_{1-n+k}=d_k,\;k=1,\ldots,n}
	\subseteq \D^\Z$. We assume $P(A)>0$, as otherwise $\nu=0$. Further define the non-normalised measure
	$\Ph \defeq P(A\cap \fdot) \in \Mp(\D^\Z)$ and denote $\mu=\Ph \circ \(\condp{\Xf}{\Xp}\)^{-1}$. Note that the
	conditional probability $\condp{\Xf}{\Xp}$ in the definition of $\mu$ is w.r.t.\ to $P$ not $\Ph$.
	Using stationarity of $P$ we obtain
	\begin{eqnarray*}
		\nu &=& P_\N\([d_1,\ldots, d_n] \cap \shift^{-n}(\fdot)\) \= 
			\intp{\Bcondp{A \cap \sset{X_\N\in \fdot}}{\Xp}} \= \intap{A}{\condp{\Xf}{\Xp}} \\
		&=& \plainint{\condp{\Xf}{\Xp}}{\Ph} \= \intmu{\id_{\Pf}} \= \|\mu\|\cdot r\Bigl(\frac{\mu}{\|\mu\|}\Bigr),
	\end{eqnarray*}
	where $\id$ is the identity, $\|\mu\|=\mu\(\Pf\)$ is the norm of total variation and $r\colon \PPf \to \Pf$ is the
	resultant (also called barycentre map) from integral representation theory. Because $\Ph\ll P$,
	and thus $\mu\ll \mueps$, the support of $\mu$ is contained in $\causupp$. Due to compactness of $\causupp$,
	this implies that the barycentre lies in the closed convex hull of $\causupp$ (\cite{ChoquetII}), i.e.\ 
		\[ r\(\tfrac{1}{\|\mu\|}\mu\) \in \wsclosure{\conv(\causupp)} \und \nu\in\wsclosure{\linhull(\causupp)}. \]
\wcase{``$\supseteq$''} We have to show that $\wsclosure{V_P}$ has full $\mueps$-measure, in other words that
	$\condp{\Xf}{\Xp} \in \wsclosure{V_P}$ $P$-a.s. By the martingale convergence theorem we have for all $B\in\B(\future)$
	a.s.:
		\[ \scondp{\Xf\in B}{\Xp}(\omega) \= \nlim \scondp{\Xf\in B}{X_{[-n,0]}}(\omega)
			\= \nlim\, \frac{\tauD_{X_1(\omega)\cdots X_n(\omega)}(P_\N)(B)}{P\([X_1(\omega),\ldots, X_n(\omega)]\)} \]
	Because $\B(\future)$ is countably generated and setwise (pointwise) convergence of a sequence of probability measures
	implies {\weaks} convergence, we obtain $\condp{\Xf}{\Xp} \in \wsclosure{\R\cdot Q_P}$ $P$-a.s.
\end{proof}

\begin{corollary}\corlabel{upperbound}
	The logarithm of the process dimension is upper bounded by the topological statistical complexity (the logarithm of the
	number of causal states).
\end{corollary}

As a second corollary, we obtain an ergodic decomposition formula for process dimension. Namely, the dimension of a process is
the sum of the dimensions of its ergodic components. This is not too surprising, because ergodic measures are mutually singular.
We prove this formula more generally in the not necessarily commutative case in \secref{results}. Nevertheless, we give an
alternative proof for the classical case here.

\begin{corollary}\corlabel{classergdecomp}
	Let $P\in\PO$ with ergodic decomposition $\nu\in\P\(\Pe\)$. Then 
		\[ \dim(P) \= \sum_{\mu\in\supp(\nu)} \dim(\mu), \]
	where we use the convention that sums over uncountably many strictly positive elements are infinite.
\end{corollary}
\begin{proof}
	We use that $\dim(P)=\dim\(\linhull(\causupp)\)$ by \propref{closeq}. It is evident that $\dim(P)$ cannot
	exceed the sum. Let $P_1,\ldots,P_n\in\supp(\nu)$ be distinct ergodic components of $P$. Then there are
	disjoint $A_1,\ldots,A_n\in\B(\future)$ s.t.\ $P_k\(\sset{\Xf\in A_k}\)=1$.  Consequently,
	$P_k\(\sset{\Xf\in A_k}\bigmid \Xp\)=1$ $P_k$-a.s.\ and, because $M_k=\set{\mu\in\Pf}{\mu(A_k)=1}$ is
	closed, $\caus_{P_k}\subseteq M_k$. The vector spaces $\linhull(M_k)$ are obviously linearly independent
	and thus the vector spaces $V_k=\linhull(\caus_{P_k})$ are linearly independent as well.
	As $\linhull(\causupp) \supseteqs \bigcup_k V_k$, we obtain $\dim(P) \ges \sum_k \dim(P_k)$.
\end{proof}

\begin{remark}
	Assume that $\causupp$ is countable and all elements have non-zero $\mueps$-probability, so that we can identify
	$\causupp$ with the set of causal states. The $\eps$-machine of computational mechanics is an HMM with the set of causal
	states as internal states. The vector space corresponding to this HMM (if we interpret it as OOM) is $\Ms(\causupp)$ as
	opposed to the canonical OOM vector space $V_P=\linhull(\causupp)$. The latter can be much lower dimensional, because it
	utilises the linear structure of $\causupp$.
\end{remark}

The closures in \propref{closeq} are really necessary, as we see in the next example.
Although $\causupp$ is closed, $\linhull(\causupp)$ is not (in general). Also, in general, neither does
$\linhull(\causupp)$ contain $V_P$ nor the other way round.

\begin{ex}\exlabel{infinitedim}
	Let $\D=\sset{0,1}$ and for $p\in[0,1]$ let $P_p\in\PO$ be the Bernoulli process with parameter $p$, i.e.\ $P_p$ is
	i.i.d.\ with $P_p\(\sset{X_1=1}\)=p$. Consider the uncountable mixture $P=\plainint{P_p}{p}$, where integration is
	w.r.t.\ Lebesgue measure. Then $\mueps$ is the image of Lebesgue measure under the map $p\mapsto P_p\circ\Xf^{-1}$,
	and $\causupp \= \bset{P_p\circ \Xf^{-1}}{p\in[0,1]}$ is the set of i.i.d.\ processes. We make the following observations:
	\begin{enumerate}[1.]
		\item $\linhull(\causupp) \cap \Pf$ is the set of \emph{finite} mixtures of i.i.d.\ processes, in particular
			$\linhull(\causupp)$ is \emph{not} closed.
		\item $V_P$ has countable algebraic dimension, i.e.\ it is the linear hull of a countable set, while a basis of 
			$\linhull(\causupp)$ has to be uncountable (the family $(P_p)_{p\in[0.1]}$ is linearly independent).
			Thus, $V_P$ cannot contain $\linhull(\causupp)$.
		\item All elements of $V_P\cap \Pf$ have an uncountable number of ergodic components. Therefore,
			$\linhull(\causupp)$ and $V_P$ are even disjoint.\exend
	\end{enumerate}
\end{ex}

\section{Non-Commutative OOMs}

Since OOMs are, unlike the concept of causal states, formulated
algebraically, they have a rather natural generalisation to the
setting of non-commutative algebras: Intuitively, we have to replace
symbols from an alphabet $\D$ by operators representing
observables. More precisely, we pass from the algebra $\cont(\D)$ of
(continuous) complex functions on Delta to an operator algebra $\A$
with their self-adjoint operators usually associated with observables
of a quantum system. In \cite{Fannes}, corresponding models, here
referred to as NC-OOMs, have been introduced and investigated in
detail for a class of stationary states on quasi-local C*-algebra that
feature finite process dimension as introduced in
\defref{NC-processdim} below. In what follows we do not impose this
restriction.

Let $\A$ be a \emph{finite-dimensional} C*-algebra with unit $\unita$
and positive cone $\A_+=\set{a\in\A}{a\ge 0}$. With $\Af\defeq
\tensor\A\N$, we denote the C*-algebraic tensor product (i.e.\ the
norm completion of the algebraic tensor product), and similarly,
$\AI\defeq \tensor\A{I}$ for $I\subseteq\Z$.  Let $\State(\A)$ denote
the set of states on $\A$, i.e.\
$\State(\A)=\bset{\rho\in\A^*}{\text{$\rho$ positive, }
\rho(\unita)=1}$, where $\A^*$ denotes the dual space of $\A$.  Note that if $\,\A=\cont(\D)=\C^\D$, then
$\Af=\cont(\future)$ and $\State(\Af)$ can be identified with $\P(\future)$.

\begin{definition}\deflabel{ncoom}
	$(V, T, v, \ell)$ is an \define{NC-OOM} with output algebra $\A$ if $V$ is a vector space,
	$T\colon \A \times V \to V$, $(a, w) \mapsto T_a(w)$ is bilinear, $v\in V$ and $\ell\in V^*$ such that
	for $n\iN$, $a_1,\ldots, a_n\in\A_+$
	\begin{gather*}
		\eqnitem{1} \ell(v) \= 1, \sepeqnitem{2} \ell\circ T_{\unita} \= \ell,
		\sepeqnitem{3} \ell\circ T_{a_n}\circ\cdots\circ T_{a_1}(v) \ges 0.
	\end{gather*}
	The state $\vphi$ on $\Af$ obtained by linear extension of $\vphi(a_1 \otimes \cdots \otimes a_n)
	\defeq \ell\circ T_{a_n}\circ\cdots\circ T_{a_1}(v)$ is called \define{generated} by the NC-OOM.
\end{definition}

\begin{enremark}
	\item The state $\vphi$ generated by an NC-OOM is a well defined state on $\Af$. Note that it is not
		necessarily translation invariant.
	\item We adapted the definition of finitely correlated states given in \cite{Fannes}, to fit the
		classical OOM definition, see \defref{oom}. In \cite{Fannes}, the $T_{a_k}$ are applied in
		reverse order and $\ell$ need not be normalised: there
		$\vphi(a_1\otimes \cdots \otimes a_n) \=
			\frac1{\ell(e)} \cdot \ell\circ T_{a_1}\circ\cdots\circ T_{a_n}(v)$.
		Note that as a consequence of the reverse order combined with condition 2.\ the associated
		finitely correlated states in \cite{Fannes} are translation invariant by construction: 
		\begin{eqnarray*} 
		\varphi(\unita\otimes a_1 \otimes \cdots \otimes a_n) &=&
			\frac{1}{\ell (v)} \ell \circ T_{\unita} \circ T_{a_1} \circ \cdots \circ T_{a_n}(v)\\
		  &=& \frac{1}{\ell (v)} \ell \circ T_{a_1} \circ \cdots \circ T_{a_n}(v)\\
		  &=& \varphi(a_1 \otimes \cdots \otimes a_n)
		\end{eqnarray*} 
\end{enremark}

Let $\MCf$ be the set of complex-valued measures of bounded variation on $\future$, and $\psi\colon\MCf \to
\cont(\future)^*$ the natural isomorphism. Then every probability measure $P\in\Pf$ corresponds to the state
$\psi(P)$ on the commutative C*-algebra $\cont(\future)$.
In the same vein, OOMs with output alphabet $\D$ can be interpreted as the special case of NC-OOMs with
commutative output algebra $\A=\cont(\D)$.
More precisely, there is a natural one-to-one correspondence $\iota$ as follows.
If $\Oom\= \(V, (T_d)_{d\in\D}, v, \ell\)$ is an OOM, the corresponding NC-OOM is
$\iota(\Oom) \= (\Vt, T, v, \ellt)$, where $\Vt\= V\oplus iV$ is the complexification of $V$, and $\ellt$ is
the complex-linear extension of $\ell$ to $\Vt$. $T$ is given by $T(f, w) \defeq \sum_{d\in\D} f(d) T_d(w)$
for $f\in\A=\cont(\D),\; w\in V$, and extended linearly to $w\in \Vt$. Obviously, $\dim_\R(V) \=
\dim_\C(\Vt)$. Furthermore, it is straight-forward to check that if $\Oom$ generates $P\in\Pf$, denoted by
$\gen(\Oom)=P$, then $\iota(\Oom)$ generates $\psi(P)\in\Sf$, denoted by $\gen\(\iota(\Oom)\) \= \psi(P)$.
This means that the following diagram commutes:
\newcommand{\bij}{\ar@{^{(}->>}}\newcommand{\sur}{\ar@{>>}}
\[\xymatrix{
	\text{OOMs}(\D) \bij[r]^-{\iota} \sur[d]^-{\gen} & \text{NC-OOMs}\(\cont(\D)\) \sur[d]^-{\gen} \\
	\Pf \bij[r]^-\psi & \State\(\cont(\future)\)
}\]

The canonical NC-OOM of a state $\vphi\in\Sf$ is defined similarly to the canonical OOM of a classical
probability distribution, cf.\ \defref{canoom}. In more detail, the dual $\Af^*$ of $\Af$ corresponds to the
space $\Msf$ of signed measures used in the classical construction. The initial vector is $\vphi$ itself, and
the evaluation functional $\ell_\A$ is the evaluation at $\unitaf$, i.e.\ $\ell_\A(\rho)=\rho(\unitaf)$. The map
$\tauA \colon \A\times\Af^* \to \Af^*$ is defined by
	\[ \tauA(a, \rho)\defeq \tauA_a(\rho) \defeq \rho(a\otimes\fdot) \defeq \( X \mapsto \rho(a\otimes X) \), \]
and again we set
	\begin{eqnarray}\eqlabel{tau-fct} \tauA_{a_1\cdots a_n} \defeq \tauA_{a_n}\circ\cdots\circ\tauA_{a_1}. \end{eqnarray}

\begin{definition}\deflabel{canNCOOM}
	For $\vphi\in\Sf$ let
		\[ V_\vphi \defeq \linhull
			\bset{\tauA_{a_1\cdots a_n}(\vphi)}{n\iN_0,\; a_1,\ldots,a_n\in\A} \subseteqs \Af^* \]
	and denote the function $\A\times V_\vphi\to V_\vphi$, $(a,\rho)\mapsto\tauA_a(\rho)$ with a slight abuse of notation
	again by $\tauA$.  Then $(V_\vphi, \tauA, \vphi, \ell_\A)$ is called \define{canonical NC-OOM} of $\vphi$.
\end{definition}

\begin{enremark}
	\item $V_\vphi$ is a vector space and $\tauA_a$, for all $a\in\A$, maps $V_\vphi$ into $V_\vphi$.
	\item The discussion in \cite{Fannes} is about translation invariant states on $\AZ$. There, the image
		$W_\vphi$ of the map $\Ap\to\Af^*$, $a\mapsto \vphi(a\otimes \fdot)$  is used instead of $V_\vphi$. In general,
		we have the relation $V_\vphi\subseteqs W_\vphi\subseteqs \wsclosure{V_\vphi}$. In the finite-dimensional case,
		however, the two spaces coincide.
	\item Let $P\in\Pf$, and $\Oom$ the canonical OOM of $P$. Then the corresponding NC-OOM $\iota(\Oom)$
		is the canonical NC-OOM of $\psi(P)$ up to the identification of $\MCf$ with $\cont(\future)^*$
		by the isomorphism $\psi$. In particular,
		$\tauA_a \= \sum_{d\in\D} {a(d)\cdot \psi\circ \tauD_d \circ \psi^{-1}}$ and
			\[ \dim(P) \= \dim\(\psi(P)\). \]
	\item Note that $\tauA_a$ is {\weaks} continuous.
\end{enremark}

\begin{lemma}
	The canonical NC-OOM of\/ $\vphi\in\Sf$ generates\/ $\vphi$.
\end{lemma}
\begin{proof}
	Let $a_1,\ldots, a_n\in\A$. We obtain
		\[ \ell_\A\circ \tauA_{a_1\ldots a_n} (\vphi) \= \tauA_{a_n}\(\tauA_{a_1\ldots a_{n-1}}(\vphi)\) (\unitaf)
			\= \tauA_{a_1\ldots a_{n-1}}(\vphi) (a_n) \= \cdots \= \vphi(a_1\otimes \cdots \otimes a_n). \qedhere \]
\end{proof}

Similarly to the definition of process dimension of a probability distribution as given in \defref{procdim}, we propose:
\begin{definition}\deflabel{NC-processdim}
	The \define{process dimension} of $\vphi\in\Sf$ is the dimension of its canonical NC-OOM:
		\[ \dim(\vphi) \defeq \dim(V_\vphi)\ins\N\cup\sset{\infty}. \]
\end{definition}

\section{Properties of Process Dimension}\seclabel{results}

In this section we present our main results: lower semi-continuity and
an ergodic decomposition formula for process dimension. In the
classical special case, corresponding results for a class of complexity measures 
have been obtained in
\cite{thesis, funcpaper}. For the technical prerequisits that are required for our 
non-commutative extension we refer to the books \cite{BratteliRobinson, Ruelle}.

\begin{theorem}
	The process dimension\/ $\dim\colon \State(\Af) \to \N\cup\sset{\infty}$ is {\weaks} lower semi-continuous.
\end{theorem}
\begin{proof}
	Because $\Af$ is separable, $\Sf$ is {\weaks} metrisable and thus sequential semi-continuity implies semi-continuity.
	Let $\vphi$ be the \weaks\ limit of a sequence $\folge{\vphi}$ in $\Sf$ and $\dim(\vphi)\ge d$. We have to show
	that $\dim(\vphi_n)\ge d$ for sufficiently large $n$.
	Let $(V_\vphi, \tauA, \vphi, \ell_\A)$ be the canonical NC-OOM of $\vphi$. Since, by \defref{canNCOOM}, 
	$\dim(V_\vphi)\=\dim(\vphi)$, we can choose linearly independent $v_1,\ldots,v_d\in V_\vphi$. Moreover, by definition
	of $V_\vphi$, there exist $a_{kj}\in\A$, $1\le k\le d$, $1\le j\le m_k$, such that
	$v_k=\tauA_{A_{k1}\ldots A_{km_k}}(\vphi)$.
	For $n\iN$, we define vectors $v^n_k \defeq\tauA_{A_{k1}\ldots A_{km_k}}(\vphi_n)$ in $V_{\vphi_n}$,
	respectively. Due to continuity of $\tauA_A$, we have $v^n_k \tows v_k$.
	If $v^n_1,\ldots,v^n_d$ are linearly independent for all sufficiently large $n$, the proof is finished.
	Suppose this is not the case and w.l.o.g.\ that they are dependent for all $n$. Then there are $\lambda^n_k\in[-1,1]$
	with $\max_k |\lambda^n_k| =1$ and $\sum_{k=1}^d \lambda^n_k v^n_k = 0$ for all $n$. Because $[-1,1]^d$ is compact, we
	may assume by passing to a subsequence that $\lambda^n_k\ton \lambda_k$ for some $\lambda_k$. Due to {\weaks} continuity
	of addition and scalar multiplication, $\sum_k \lambda_k v_k=0$ and hence $\lambda_k=0$ for all $k$, in contradiction to
	$\max_k |\lambda^n_k| = 1$.
\end{proof}

Due to the one-to-one correspondence between canonical OOMs of classical processes and canonical NC-OOMs of
associated states on (abelian) C*-algebras the above theorem has the following corollary.
\begin{corollary}
	The classical process dimension\/ $\dim\colon\Pf\to \N\cup\sset{\infty}$ is {\weaks} lower semi-continuous.
\end{corollary}

We now derive an ergodic decomposition formula for process dimension
in the stationary case. Let $\Ssf$ be the convex set of translation
invariant states, and $\Sef\subseteq\Ssf$ the set of ergodic states,
i.e.\ extreme points in $\Ssf$. Since $\Af$ is asymptotically abelian
(w.r.t.\ the shift) and $\Ssf$ is metrisable, $\Ssf$ is a simplex and
the set $\Sef$ of ergodic states is measurable in $\Ssf$. In
particular, every translation invariant state $\vphi$ has a unique
ergodic decomposition $\nu\in\P\(\Ssf\)$, which is supported by the
ergodic states, $\nu(\Sef)=1$, and $\vphi$ is the barycentre of $\nu$,
	\[ \vphi \= \intanu{\Sef}{\id}. \]
Moreover, in what follows, we make use of the important fact that in our situation the ergodic decomposition
is orthogonal. For details, see \cite[Sec.~4.1, 4.3.1]{BratteliRobinson}. We obtain the following ergodic
decomposition formula for process dimension.

\begin{theorem}\thmlabel{ergdecomp}
	Let\/ $\vphi\in\Ssf$ be a translation invariant state with ergodic decomposition\/ $\nu\in\P\(\Sef\)$. Then
		\[ \dim(\vphi) \= \sum_{\psi\in\supp(\nu)} \dim(\psi). \]
\end{theorem}

For the proof, we use the following two lemmas.

\begin{lemma}[representation on Hilbert space]\lemlabel{representation}
	Let\/ $\vphi\in\Ssf$ and\/ $(\Hphi, \piphi, \xiphi)$ be the GNS-representation of\/ $\Af$ w.r.t.\ $\vphi$.
	Then there is a linear injection\/ $\rep$ from\/ $V_\vphi$ into\/ $\Hphi$ with
		\[ \rho(X) \= \bscalphi{\piphi(X)\rep(\rho)}{\xiphi} \for X\in\Af,\, \rho\in V_\vphi \]
\end{lemma}
\begin{proof}	
	Let $u$ be the unitary representation of the shift on $\Af$. Then stationarity of $\vphi$ implies $u\xiphi=\xiphi$. For $a_1,\ldots,a_n\in\A$ and $A=a_1\otimes \cdots \otimes
	a_n$ let
		\[ \xi_A \defeq u\piphi(a_n)\cdots u\piphi(a_1)\xiphi \= u^n \piphi(A)\xiphi \ins \Hphi\]
	and extend the definition to $A\in \AI[\sset{1,\ldots,n}]$ linearly.
	Then we have for $A=\sum_{i=1}^m a_{i1}\otimes\cdots\otimes a_{in}$
	\begin{equation}\eqlabel{taurep}
		\tauA_A(\vphi)(X) \defeq \sum \tauA_{a_{i1}\cdots a_{in}}(\vphi)(X)
			\=  \bscalphi{(u^n)^*\piphi(X)u^n\piphi(A)\xiphi}{\xiphi}
			\EQ{$u\xiphi=\xiphi$} \bscalphi{\piphi(X) \xi_A}{\xiphi},
	\end{equation}
	where $\tauA_{a_{i1}\cdots a_{in}}$ is defined in \eqref{tau-fct}.
	For every $\rho\in V_\vphi$, there is an $n\iN$ and $A\in\Afin\defeq \bigcup_{n\iN}
	\AI[\sset{1,\ldots,n}]$ with $\rho\=\tauA_A(\vphi)$. Define $\rep(\rho)\defeq \xi_A$.
	Because $\xiphi$ is cyclic and $\Afin$ is dense in $\Af$,
	$\bscalphi{\piphi(X)\zeta_1}{\xiphi}\=\bscalphi{\piphi(X)\zeta_2}{\xiphi}$ for all $X\in\Afin$ implies
	that $\zeta_1=\zeta_2$. Thus $\rep$ is well-defined. Injectivity is obvious, because $\rho$ can be
	recovered from $\xi_A$ by \eqref{taurep}.
\end{proof}

\begin{lemma}\lemlabel{sum}
	Let\/ $\vphi\in\Ssf$ and\/ $\vphi=\sum_{\psi\in \Psi} \nu(\psi) \psi$, where $\Psi\subset\Ssf$ is countable and
	$\nu(\psi) \gs 0$.  Then
		\[ \wsclosure{V_\vphi} \= \wsclosure{\sum_{\psi\in \Psi}\, V_\psi}, \]
where $V_\vphi$ is defined in \defref{canNCOOM}.
\end{lemma}
\begin{proof}
\wcase{``$\subseteq$''} By linearity of $\tauA_a$, we obviously have $V_\vphi \subseteq \sum_\psi V_\psi$.
\wcase{``$\supseteq$''} Because $\nu(\psi)\psi \les \vphi$, there is a $\xi_\psi\in\Hphi$ with $u\xi_\psi=\xi_\psi$ and
	$\psi(X) \= \bscalphi{\piphi(X)\xi_\psi}{\xiphi}$. As $\xiphi$ is cyclic and $\Afin$ is dense in $\Af$, there is
	a sequence $A_n\in\AI[\sset{1,\ldots,n}]$ with $\piphi(A_n)\xiphi \to \xi_\psi$. Let $\rho_n\defeq
	\tauA_{A_n}(\vphi)$. Then $\rho_n\in V_\vphi$, and for all $X\in\Af$
	\[
		\rho_n(X) \= \bscalphi{\;\underbrace{\piphi(X)u^n}_{\mathclap{\|\cdot\| \= \|\piphi(X)\|\ls\infty}}\;
				\piphi(A_n)\xiphi}{\xiphi}
			\tons \bscalphi{\piphi(X)\xi_\psi}{\xiphi} \= \psi(X).
	\]
	Thus $\rho_n\tows \psi$ and $\psi\in\wsclosure{V_\vphi}$. As $\tauA_a$ is \weaks\ continuous,
	$\tauA_a\(\wsclosure{V_\vphi}\) \subseteq \wsclosure{V_\vphi}$, hence  $V_\psi \subseteq \wsclosure{V_\vphi}$.
\end{proof}

The lemma shows that in order to represent $\vphi$ in terms of NC-OOMs, we have to represent all ergodic components $\psi$ of
$\vphi$ (and not more). We still have to show that they have to be represented independently without synergies. This follows
easily from the orthogonality of the ergodic decomposition.

\begin{proof}[Proof of \thmref{ergdecomp}.]\begin{enumerate}
\itcase{Finitely many ergodic components}
	In this case, \lemref{sum} directly implies ``$\le$.'' For ``$\ge$,'' we may assume that $\dim(V_\vphi)\ls \infty$ and
	thus also $\dim(V_\psi)\ls\infty$ for all ergodic components $\psi$. In particular, $V_\vphi \= \sum_\psi V_\psi$ by
	\lemref{sum}.
	We can identify the GNS-Hilbert space $\H_\psi$ with a subspace of $\Hphi$ and because the ergodic decomposition is
	orthogonal, the $\H_\psi$ are mutually orthogonal. Since $\rep(V_\psi) \subseteq \H_\psi$, by
	\lemref{representation}, the sum of vector spaces is direct, i.e.\ $V_\vphi \= \bigoplus_\psi V_\psi$,
	and $\dim(V_\vphi) \= \sum_\psi \dim(V_\psi)$.
\itcase{Infinitely many ergodic components}
	The sum on the right-hand side is infinite. To see that also $\dim(\vphi)\=\infty$, fix $n\iN$ and choose a
	decomposition of $\supp(\nu)$ into disjoint measurable subsets $\Psi_k$, $k=1,\ldots,n$ with positive $\nu$-measure.
	Define $\psi_k\defeq \intanu{\Psi_k}{\id}$ to be the barycentre of $\nu\restricted{\Psi_k}$. Then the decomposition
	$\psi\=\sum_k \psi_k$ is orthogonal and by the above argument $\dim(V_\vphi) \= \sum_k \dim(V_{\psi_k}) \ges n$.\qedhere
\end{enumerate}\end{proof}

Note that \thmref{ergdecomp} provides an alternative proof of \corref{classergdecomp}.

\begin{acknowledgement}
	This work has been supported by the Santa Fe Institute.
\end{acknowledgement}

\bibliography{bib-math}

\end{document}